\documentclass[a4paper,12pt]{article}
\usepackage[toc,page]{appendix}
\usepackage{fancyhdr,lastpage}
\usepackage{amsthm}
\usepackage{amsmath}
\usepackage{amssymb}
\usepackage{color}
\numberwithin{equation}{section}
\newtheorem{theor}{Theorem}[section]
\newtheorem{prop}{Proposition}[section]
\newtheorem{lem}{Lemma}[section]

\newtheorem{cor}{Corollary}[section]
\newtheorem{rem}{Remark}[section]
\newcommand{\eps}{\varepsilon}
\newenvironment{Proof}[1][Proof.]{\begin{trivlist}
\item[\hskip \labelsep {\bfseries #1}]}{\end{trivlist} }

\numberwithin{equation}{section}

\usepackage{graphics}
\usepackage{fullpage}
\usepackage{url}
\usepackage{enumerate}
\usepackage{hyperref}
\usepackage{pgf}
\title{Boundedness  of global solutions of a $p$-Laplacian evolution equation with a nonlinear gradient term}
\author{Amal Attouchi}
\date{\footnotesize\textsl{Universit\'e Paris 13, Sorbonne Paris Cit\'e, Laboratoire Analyse, G\'eom\'etrie et Applications, CNRS, UMR 7539, 93430 Villetaneuse, France}}

\begin{document}

\maketitle
\begin{abstract}
We investigate the boundedness and large time behavior
of  solutions  of the Cauchy-Dirichlet problem for the  one-dimensional degenerate parabolic equation with gradient nonlinearity:
$$u_t=(|u_x|^{p-2} u_x)_x+|u_x|^q\qquad \text{in}\quad (0,\infty)\times(0,1), \qquad q>p>2.$$
We prove that: either $u_x$ blows up in finite time, or $u$ is global and converges in $W^{1,\infty}(0,1)$  to the unique steady state. This in particular eliminates the possibility of global solutions with unbounded gradient. For that purpose a Lyapunov functional is constructed by the approach of Zelenyak. 

\end{abstract}

\section{Introduction and main results}
In this paper we are interested in the asymptotic behavior of global solutions to the following  one-dimensional degenerate diffusive Hamilton-Jacobi equation
\begin{equation}\label{eq1}
 \left\{
\begin{array}{lll}
u_t-(|u_x|^{p-2} u_x)_x=|u_x|^q, &0<x<1,  t>0,\\
u(t,0)=0,\qquad u(t,1)=M\geq 0, &\qquad\qquad\quad t>0,\, \\
u(0,x)=u_0(x), &0<x<1,
\end{array}
\right.
\end{equation}
with $q>p>2$, $M\geq 0$ and suitably regular initial data $u_0$.

Problem \eqref{eq1} models a variety of physical phenomena which arise for example in the study
 of surface growth where  a stochastic version of it is known as the Kardar-Parisi-Zhang equation $(p=2, q=2)$. It has also a mathematical  interest through the viscosity approximation of Hamilton-Jacobi type equations from control theory.

Solutions of \eqref{eq1} exhibit a rich variety of qualitative behaviors, according to the values of $p\geq 2$ and  $q\in (0,\infty)$.

If $q\leq p$, it is known that all solutions are global and bounded in $W^{1,\infty}$ norm \cite{lady}.  For $q\in [p-1,p]$ it was proved in \cite{oups} that nonnegative viscosity solutions of \eqref{eq1} with homogeneous Dirichlet boundary condition decay to 0 and the rate of convergence was also obtained, see also \cite{decay} for the semilinear case. Concerning the large time behavior of global weak solutions to \eqref{eq1}  with homogeneous boundary conditions and $q\in (0, p-1)$, it has been shown that there exists a one parameter family of nonnegative steady states, and any solution converges uniformly to one of these stationary solutions (cf. \cite{stin,barle, lauren2}).

For $q>p\geq 2$, the situation is quite different. It is known that for any $M\geq 0$ and suitably large $u_0$, there exist solutions of \eqref{eq1} for which the $L^{\infty}$ norm of the gradient blows up in finite time (the $L^{\infty}$ norm of the solution remaining bounded) \cite{souma, amal}, while there exist global  and decaying solutions for  $u_0$ sufficiently small \cite{zhan}. In view of
a classification of all solutions of \eqref{eq1}, it is then a natural question to ask
whether or not $C^1$-unbounded global solutions may exist. The question of the boundedness of global solutions of \eqref{eq1} was initiated for the semilinear case $p=2$ in
\cite{Arri} and further investigated in \cite{vaz,zhan}. Denoting $M_c:= (q-1)^{\frac{q-2}{q-1}} / (q-2)$, the result of \cite{Arri} says that if $0\leq M< M_c$, then any global solution of \eqref{eq1} is bounded in $ C^1$ norm for $t \geq 0$, that is,
\begin{equation}
 \underset{t\geq 0}{\text{sup}}\quad |u_x(t,.)|_{\infty}< \infty.
\end{equation}
On the other hand, it is known from \cite{vaz} that some unbounded global solutions do
exist if $M=M_c$ and $u_0\leq  U(x)$ where $U(x):=  M_c x^{(q-2)/(q-1)}$  is the unique singular steady state. Moreover the precise exponential rates of  the gradient blow-up in infinite time was obtained.

Motivated by the results of the papers \cite{Arri, vaz}, we modify the method used by Arrieta, Rodriguez-Bernal and Souplet and extend their results  on the  classification of large time behavior of global solutions to the degenerate parabolic equation case $p>2$. 
\newline
From now on, we assume that $q>p>2$. By a solution of \eqref{eq1}, we mean a weak solution (see Section 2 below for a precise definition and well-posedness results). We recall that weak solutions of \eqref{eq1} satisfy a comparison principle, hence in particular  
\begin{equation}\label{compaprin}
\min \left\{ \underset{[0, 1]}{\text{min}}\, u_0, 0\right\}\leq u(t,x)\leq\max \left\{\underset{[0, 1]}{\text{max}}\, u_0, M\right\}, \qquad 0\leq t< T_{max}, \qquad 0\leq x\leq 1.
\end{equation} 
Our main result is then the following:

\begin{theor}\label{Main result}
Assume that  $q>p>2$ and $u_0$ in $W^{1, \infty}(0,1)$, $u_0(0)=0, u_0(1)=M$. Set $M_b=\frac{q-p+1}{q-p} \left(\frac{q-p+1}{p-1}\right)^{\frac{1}{p-1-q}}$.
\begin{enumerate}[(i)]
\item If $0\leq M<M_b$, then any global weak solution  of \eqref{eq1} is bounded in $C^{1}$ norm. Moreover it converges in $C^1([0,1])$ to the unique steady state.
\item If $M>M_b$, then all  weak solutions of \eqref{eq1} exhibit gradient blow-up in finite time.
\end{enumerate}
\end{theor}

The proof of Theorem \ref{Main result}  proceeds by contradiction. It relies on the analysis of steady states and the existence of a Lyapunov functional which enjoys nice properties on any global trajectory of \eqref{eq1}, even if it were unbounded in $C^{1}$ norm.
The construction of such a nice Lyapunov functional which is handled through the Zelenyak technique, together with the fact that the singularities may only take place near the boundary, allows  us to prove the following convergence result: any global solution, even unbounded in $W^{1,\infty}$ must converge in $C([0,1])$ to a stationary solution $W$ of \eqref{eq1} with $W(0)=0$, $W(1)=M$ (see Proposition \ref{converg}). On the other hand, if $u$ were unbounded, then our gradient estimates would imply that $W_x(0)=+\infty$ or $W_x(1)=-\infty$. But such a $W$ does not exist if $M\neq M_b$, leading to a contradiction.

Although the scheme of proof follows that in \cite{Arri} for $p=2$, we have to face a number of additional technical difficulties, caused by the lack of regularity of solutions. In particular, we have to work at the level of regularized problems, including for the construction of Lyapunov-Zelenyak functional. This, in turn requires good convergence properties and estimates of regularized solutions. For this, we heavily rely on results from our previous work \cite{amal} (which concerned the higher dimensional problem as well) and an extension up to the boundary of a result of DiBenedetto-Friedman on  the regularity  of the derivative of weak solutions of a degenerate parabolic problem (see Proposition \ref{friedmanetendu}).

\begin{rem}
\begin{enumerate}[(a)]

\item For the critical case $M=M_b$, all solutions must blow up  in either finite or infinite time.  The existence of global solution which are unbounded in $W^{1,\infty}$ norm (that is infinite time gradient blow-up) should occur for some suitable initial data as it is the case for the corresponding semilinear equation (namely for initial data $u_0$
below the singular steady state), but this still is an open problem. Moreover,  we know  from Proposition \ref{converg} that, even in this case the solutions will converge in $C([0,1])\cap C^1_{loc}((0,1])$ to the unique singular steady state.

\item Since the technique of Zelenyak to obtain a Lyapunov functional is restricted to the one-dimensional setting, the large time behavior and  the boundedness of global solution in $W^{1,\infty}$ norm are still open problems in higher dimension.
\end{enumerate}

\end{rem}

Let us mention some results  concerning related equations possessing solutions with unbounded gradient. When the nonlinearity is replaced with an exponential one and $p=2$, results on boundedness and existence of infinite time gradient blow-up solutions are obtained in \cite{expo1, expo2}.  A phenomenon of infinite time gradient blow-up has been observed  for quasilinear equations involving mean curvature type operators \cite{chen1997}. For results on  interior gradient blow-up we refer the reader to \cite{ang,asai}. Finally for other results concerning existence,  asymptotic behavior of global solutions for the corresponding Cauchy problem and a viscosity solution approach see \cite{chen, lauren1, superlinear} and references therein.
\bigskip

The rest of the paper is organized as follows. Section 2 contains some useful
preliminary material, including smoothing properties of
solutions and estimate of the derivative $u_x$. In section 3, we employ the technique of Zelenyak \cite{zelen}, along with a trick used in \cite{Arri}, to construct an approximate Lyapunov functional for weak solutions to \eqref{eq1}. Section 4 is devoted to the proof of Theorem \ref{Main result}.

\section{Preliminary  estimates and steady states}
\subsection{Space and time derivative estimates}

For $u_0\in W^{1,\infty}((0,1))$, $u_0(0)=0, u_0(1)=M$, by a (weak) solution of \eqref{eq1}  on $[0, T]$, we mean
a function $u\in C([0,T)\times[0,1])\cap L^q((0,T); W^{1,q}(0,1))$ such that
$$ u_t\in L^2((0,T); L^2(0,1)),\quad u(0, x)=u_0(x),\quad  u(t,0)=0,  u(t,1)=M $$ and

\begin{equation}\label{defmil}
 \int_0^T\int_0^1 u_t \psi+|u_x|^{p-2} u_x\cdot \psi_x \, dx\, dt=\int_0^T\int _0^1 |u_x|^q \psi \,dx \,dt,
\end{equation}
holds for all  $\psi\in C^0([0, T]\times[0,1]) \cap L^p((0,T); W^{1,p}_0((0,1))$.

\begin{rem}
if $u$ solves \eqref{eq1}, then so does $v$ defined
by $v(t, x) = M -u(t,1-x)$.
\end{rem}
It is known (see e.g., \cite{amal}) that there exists $T_{max}= T_{max}(u_0) \in (0,\infty]$ such that
for each $T \in(0, T_{max})$, \eqref{eq1}  admits a unique solution $u$ such that  $u\in L^{\infty}\big( (0,T); W^{1,\infty}(0,1)\big)$.  In the rest of this paper, the maximal weak solution of problem \eqref{eq1} will refer to this solution.\\

Now let us state the following result (which will be very useful in the sequel) on the H\"older regularity of the derivative of solutions to a possibly degenerate parabolic problem. This result is  an extension up to the boundary (in one space dimension) of an interior estimate of DiBenedetto-Friedman \cite{dib}  (see the appendix for a proof).
The definition of a weak solution of \eqref{friedequa} is the same as in \eqref{defmil},
with $|u_x|^q$ replaced by $F$ and $u\in L^q(0,T;W^{1,q}(0,1))$ replaced by $v\in L^p(0,T;W^{1,p}(0,1))$.
\begin{prop}\label{friedmanetendu}

Let $\eps\in[0,1), u_0\in W^{1,\infty}$, $C>0$ and $F\in L^{r} ((0,T)\times(0,1))$  for some $r>2$ with $\left\|F\right\|_{L^{r}((0,T)\times(0,1))}\leq C$. Let $v$ be a weak solution of
\begin{equation}\label{friedequa}
\left\{
\begin{array}{lll}

 v_t=\left(\left(|v_x|^2 +\varepsilon^2 \right)^{\frac{p-2}{2}}v_{x}\right)_x +F(t,x),
 & t>0, & x\in (0,1),\\
v(t,0)=0, \qquad v(t,1)=M,  &t>0,& \\
v(0,x)=u_{0}(x), &  & x\in(0,1).
       \end{array}
\right.
 \end{equation}
Then, for each $\eta>0$, $v_x\in C^{\alpha}([\eta, T-\eta]\times[0,1])$ where $\alpha>0$ and the H\"older norm of $v_x$ depend only on $C, \left\|v_x\right\|_{L^p}$ and  $\left\|v\right\|_{L^{\infty}_t,L^2_x}$.
\end{prop}
Proposition \ref{friedmanetendu} can be applied to $u$ since it is bounded in $L^1(0; T;W^{1, \infty}(0; 1))$.
As a direct consequence, we get that $u$ is a $C^1$-function w.r.t. the space variable in $(0,T)\times [0, 1]$ and that its derivative $u_x$ is locally H\"older continuous.

In order to describe the asymptotic behavior, we need to collect some preliminary estimates.
  We first give the following theorem which is of independent interest. It  gives  a useful regularizing property for local solutions of \eqref{eq1} as well as  a uniform bound on $u_t$ away from $t= 0 $ for any space dimension.

Let $\Omega\subset \mathbb{R}^N$ be  a bounded domain of class $C^{2+\alpha}$ for some 
$\alpha>0$.  Consider  the following problem

\begin{equation}\label{eqN}
 \left\{
\begin{array}{lll}
 u_t- \mathrm{div}( |\nabla u|^{p-2} \nabla u)=|\nabla u|^q,  &  x\in\Omega, t>0,\\
u(t,x)=g(x),                                                  & x\in\partial\Omega, t>0,\\
u(0,x)=u_0(x),                                             & x\in\Omega,

\end{array}
\right.
\end{equation}
where the boundary data $g$ is the trace on $ \partial\Omega$  of a regular function  in $ C^2(\overline{\Omega})$,  also denoted by $g$,  and  the initial data  $u_0$  satisfies
\begin{equation}\label{initial}
u_0\in W^{1,\infty}(\Omega), \quad u_0(x)=g(x) \qquad\text{for}\quad x\in\partial\Omega.
\end{equation}

\begin{theor}\label{reg}
 Assume that $q >  p-1$ and let  $u$ be a maximal  weak solution of problem (\ref{eqN}). We have the following statements.
\begin{enumerate}[(i)]
 \item  Let  $osc(u_0) =\underset{\overline\Omega}{\max}\, u_0 - \underset{\overline\Omega}{\min}\, u_0$, then
\begin{equation}\label{derivtsup}
u_t\leq \dfrac{1}{p-2} \dfrac{ osc(u_0)}{t}
 \qquad \mathrm{in}\quad\mathcal{D}'( (0, T_{max})\times\Omega).
\end{equation}
\item   Fix $t_0\in (0, T_{max})$,  then 
\begin{equation}\label{derivtinf}
u_t\geq -\left(\dfrac{q-p+1}{p-2}\right)\underset{[0,t_0]\times\Omega}{\sup}|\nabla u|^q - \left(\dfrac{1}{p-2} \right)\dfrac{osc(u_0)}{t}\qquad\text{in}\quad \mathcal{D}'\left((0,t_0)\times\Omega\right).
\end{equation}
\item Fix $t_0\in (0, T_{max})$, then there exists $C_1=C_1\left(p,q, t_0, \underset{[0,t_0]\times\Omega}{\sup}|\nabla u|, osc(u_0)\right)>0$ such that
\begin{equation}\label{estimut}
 |u_t|\leq C_1\qquad\text{in}\quad \mathcal{D}'\left((t_0, T_{max})\times \Omega\right).
\end{equation}
\end{enumerate}
\end{theor}
\begin{proof}
$\,$\\
The initial data being bounded and the sought-for estimate being invariant by addition of a constant, we may replace $u$ with $u-B$, where $B= \min_{\overline\Omega} u_0\geq \min_{\partial\Omega} g$. Then $u\geq 0$ by the maximum principle.
(i) This  has been proved in \cite[Theorem 1.3]{amal}.\\
(ii)  Fix $t_0\in(0, T_{max})$  and let    $D=\underset{ [0, t_0]\times \Omega}{\sup} \,| \nabla u|$.  Set $$w=u_{\lambda}:=\lambda^{\gamma} u(\lambda t, x) + t L_{\lambda} \qquad1<\lambda<\dfrac{T_{max}}{t_0}, \gamma=\frac{1}{p-2}.$$ 
where  $$L_{\lambda}=(1-\lambda^{\gamma(p-1-q)}) (\lambda^{\gamma} D)^q>0.$$
Since $u\geq  0$, $\lambda \geq 1$ and $\gamma\geq 0$, we have  $w \geq u$ on $\left\{\left\{0\right\}\times\overline{\Omega}\right\} \cup\left\{(0, t_0/\lambda)\times\partial\Omega\right\}$.
On the other hand, we have on  $(0,t_0/\lambda)\times \Omega$:
\begin{equation*}
 w_t -\Delta_p w-|\nabla w|^q 
= (\lambda^{\gamma(p-1-q)} - 1)|\nabla w|^q +L_{\lambda} 
= (1-\lambda^{\gamma(p-1-q)}) [ (\lambda^\gamma D)^q - |\nabla w|^q]
\geq 0
\end{equation*}
Hence, by the comparison principle, we get that $w \geq u$ on  $(0,t_0/\lambda)\times \Omega$, that is
\begin{equation}\label{serk}
 \lambda^{\gamma} u(\lambda t,x)-u(t,x)\geq -(\lambda^{\gamma(p-1-q)}-1) (\lambda^{\gamma} D)^q
\end{equation}
Dividing (\ref{serk}) by $(\lambda-1)$ and letting $\lambda\to1^+$, we get 
$$\gamma u+ t \, u_t\geq - (q-p+1)\gamma D^q\,t \qquad\text{in} \quad \mathcal{D}'\left((0, t_0)\times\Omega\right).$$
The estimate \eqref{derivtinf} follows.\\
(iii)\, Fix $t_0\in(0, T_{max})$. By \eqref{derivtsup}-\eqref{derivtinf}, for $h>0$ small,  we have
 \begin{equation*}
 \left\|u(t_0/2+h)-u(t_0/2)\right\|_{\infty}\leq C_1 h
 \end{equation*}
where $C_1=\left( \left(\dfrac{q-p+1}{p-2}\right)\underset{(0,t_0)\times\Omega}{\sup}\,|\nabla u|^q + \left(\dfrac{2}{p-2} \right)\dfrac{osc(u_0)}{t_0}\right)>0$.\\
Due to the translation invariance of \eqref{eq1}, for $t> t_0/2$, $u(t+h, x)$ is still a solution of \eqref{eqN} with initial condition $u(h)$. Applying a comparison principle, we obtain that
$$\left\|u(t+h)-u(t)\right\|_{\infty}\leq C_1 h.$$
Since $h$ is arbitrary small, we conclude that 
 \begin{equation*}
|u_t|\leq C_1\qquad\text{in} \quad \mathcal{D}'\left((t_0/2, T_{max})\times\Omega\right).
 \end{equation*}
\end{proof}
Thanks to the upper bound of $u_t$, we derive the following
lemma giving lower and upper bounds on $u_x$, showing that $u_x$
remains bounded away from the boundary.
\begin{lem}\label{as}
 Let $u$ be a maximal weak solution of (\ref{eq1}). For all $t_0\in(0,T_{max})$, there exists $C_2=C_2(t_0, p, osc(u_0))>0$  such that for all $t\in [t_0,T_{max})$ and $0< x< 1$ ,
\begin{center}
\begin{equation}\label{fer}
  u_x(t, x)\leq\ \left(\left(\frac{q-p+1}{p-1} x\right)^{\frac{1-p}{q-p+1}}+ C_2 x\right)^{\dfrac{1}{p-1}},
\end{equation}
\end{center}

\begin{center}
\begin{equation}\label{erf}
 u_x(t,1-x) \geq - \left(\left(\frac{q-p+1}{p-1} x\right)^{\frac{1-p}{q-p+1}}+ C_2 x\right)^{\dfrac{1}{p-1}}.	
 \end{equation}
\end{center}
\end{lem}
\begin{Proof} Fix $t\in [t_0,T_{max})$ and let $y(x)=\left(|u_x|^{p-2}u_x(t,x)-C_2 x\right)^+$, 
   where $C_2= \dfrac{osc(u_0)}{(p-2)t_0}$. On any interval $(a,b)$ with $0<a<b<1$ where $y>0$, the function $y$ satisfies in the classical sense
$y'+y^{\frac{q}{p-1}}\leq 0$. Indeed, for each  $x\in (a,b)$, we have   $|u_x|^{p-2}u_x>C_2 x> C_2 a>0$ and the function $u$ is smooth at such points since the equation is uniformly parabolic \cite{lady}. Using  theorem \ref{reg} (i), we  get that  

\begin{equation}\label{correct}
y'+y^{\frac{q}{p-1}}\leq\left((|u_x|^{p-2}u_x)_x-C_2\right)+|u_x|^q\leq 0.
\end{equation} 
This implies that $y'<0$ on $(a,b)$ so that necessarily $a=0$.
Integrating inequality \eqref{correct}, it follows that $y(x)\leq \left(\frac{q-p+1}{p-1} x\right)^{\frac{1-p}{q-p+1}}$ on $( 0,b)$  and  $y(0)>0$. 
If $y\not\equiv 0$, then we can find $c=c(t)\in (0,1]$  such that $y>0$ in $(0,c)$ and $y=0$ in $[c, 1)$. Therefore we get $y(x)\leq \left(\frac{q-p+1}{p-1} x\right)^{\frac{1-p}{q-p+1}}$ on $(0, 1)$ and \eqref{fer} is readily deduced.
In the same manner, considering  $y(x)=\big(-|u_x|^{p-2}u_x(t,1-x)-C_2 x\big)^+$, we get \eqref{erf}.
\end{Proof}

\begin{rem}
Similar gradient estimates in any space dimension are already obtained in \cite{amal} using  a more technical Bernstein type argument. 
\end{rem}

The following corollary is  a direct consequence of Lemma \ref{as} that states that, when gradient blow-up occurs on the boundary, it can  only be towards $+\infty$ at $x=0$ or towards $-\infty$ at  $x=1$.
\begin{cor}
Let $u$ be a maximal weak solution of (\ref{eq1}) and $t_0\in (0, T_{max})$. There exists $C_3=C_3(t_0, p, q, osc(u_0))>0$  such that for all $t\in [t_0,T_{max})$,

\begin{equation}
 u_x(t, 0)\geq -C_3 \qquad\text{and}\qquad u_x(t,1)\leq C_3.
\end{equation}
\end{cor}

\subsection{Steady states}
It is a well-known fact that the large-time behavior of
evolution equations is closely connected to the existence and properties of the stationary states.
In  this part we are looking for nonnegative stationary solutions $W$ of (\ref{eq1}), that is  weak solution of 
\begin{eqnarray}\label{stat}
\left\{
\begin{array}{ll}
 \left(|W_x|^{p-2} W_x\right)_x+ |W_x|^q=0,   & x\in (0,1),\\
  W(0)=0,  \qquad  W(1)=M\geq 0.
  \end{array}
  \right.
\end{eqnarray}
More precisely,  $W\in C([0,1])\cap C^1(0,1)$ is  a weak solution of \eqref{stat}
  if $W(0)=0, W(1)=M$ and $W$ satisfies
\begin{equation}
 \int_0^1\left[\left(|W_x|^{p-2}W_x\right)\,\phi_x-|W_x|^q \,\phi\right] \,dx =0\qquad\text{for any}\quad\ \phi\in C^1_c(0,1).
\end{equation}
It is not difficult to show that any weak solution in the above sense is actually a classical $C^2$ solution in $(0,1)$ (for any $x_0\in (0,1)$, consider separately the cases $W_x(x_0)\neq 0$ and $W_x(x_0)=0$).
For small values of $ M\geq 0$, problem \eqref{stat} admits a
unique  weak solution $W_M = W_M(x) \in C^2([0, 1])$. Namely, this happens for $0\leq  M < M_b$, where
$M_b$ is the critical value, 
$$M_b=\frac{q-p+1}{q-p} \left( \frac{q-p+1}{p-1}\right)^{-1/(q-p+1)}.$$
More precisely,  for $M=0$ we have $W_0=0$ and for $0\leq M< M_b$,  there exists $k=k(M)\in [0,\infty)$  such that $W_M= M_b\left[\left(x+k\right)^{\frac{q-p}{q-p+1}} -k^{\frac{q-p}{q-p+1}}\right]$.
On the other hand, there is no steady state if $M > M_b$. In the critical case $M =M_b$, there still exists
a steady state $W_{M_b} = U$, given by the explicit formula $U(x)=M_b\, x^{\frac{q-p}{q-p+1}}$. $U$ belongs to $C([0, 1]) \cap C^2((0, 1])$, but it is singular in the sense that it has infinite derivative
on the left-hand boundary, $U_x(0) = \infty$.

\section{Lyapunov functional and convergence to steady states}

Since \eqref{eq1} is a degenerate problem, we do not have sufficient regularity properties of the trajectories to construct  a good smooth Lyapunov functional (which exists for one-dimensional  uniformly parabolic equations).  Hence  we first consider a regularized problem, then the main estimate which plays a key role in the proof of the convergence to steady states will be  proved by passing to the limit $\varepsilon\to 0$ in the regularizing parameter.
\subsection{Approximate problem}
Let $\varepsilon\in (0,1/2)$.  We consider the following approximate problems:
\begin{equation}\label{eq1a}
\left\{
\begin{array}{llll}
 (u_{\varepsilon})_t=\left(\left(|(u_{\varepsilon})_x|^2 +\varepsilon^2 \right)^{\frac{p-2}{2}}(u_{\varepsilon})_{x}\right)_x
+B_{\eps}( (u_{\varepsilon})_x)
 \qquad\qquad (t,x)\in(0,+\infty)\times(0,1),\\
u_{\varepsilon}(t,0)=0,  \qquad u_{\varepsilon}(t,1)=M,     \qquad t>0, \\
u_{\varepsilon}(0,x)=u_{0}(x),\qquad x\in(0,1),
       \end{array}
\right.
 \end{equation}
where $B_{\eps}(v)=\left(|v|^2+\varepsilon^2\right)^{\frac{q-2}{2}}(|v|^2+\frac{\eps^2}{p-1})$.\\
Here  we collect some useful properties of the sequence $\left\{u_{\varepsilon}\right\}$ which we will use later on.

Let  $u\in L^{\infty}\big([0,T);  W^{1,\infty}((0,1))\big)$  for  any $T\in (0, T_{max})$ be the unique, maximal  weak solution of problem (\ref{eq1}) and let $u_{\eps}$ be the unique, maximal classical solution of \eqref{eq1a} and $T(u_{\eps})$ be its existence time. We have the following proposition.
\begin{prop}\label{estimsuite}
Let $A>0$ and assume that $\left\|u_0\right\|_{W^{1, \infty}}\leq A$. Then for $0<T<T_{\max}$ and $\eps$ small, we have
\begin{enumerate}[a)]
\item $T(u_{\varepsilon})> T$,  $u_{\varepsilon} \to u$ in $C([0, T]\times [0,1])$ and $(u_{\varepsilon})_x\to u_x$ in $C_{loc}((0,T]\times[0,1])$.
\item $|(u_{\varepsilon})_x(t,x)|\leq C= C(A,T)$ on $[0, T]\times[0,1]$ and $$ \int_0^{T}\int_0^1 (u_{\varepsilon})_t^2 \ dx dt  \leq \tilde{C}(A,T).$$
\end{enumerate}
\end{prop}
\begin{proof}
For the convenience of the proof, we shall actually replace the initial data $u_0$ in the approximate problem \eqref{eq1a} with a sequence $u_{\eps,0} \in W^{1, \infty}((0,1))$, where $u_{\eps,0}\to u_0$ in $W^{1,\infty}((0,1))$,
and prove that Proposition \ref{estimsuite} remains true in this more general situation.
We know from  \cite[section 3]{amal} that there exist a small time $\tilde{\tau}=\tilde{\tau}(A)>0$  and  a subsequence  $\left\{u_{\eps_n}\right\}$ of $\left\{u_{\eps}\right\}$  such that $u_{\eps_n}$ converges in $C([0, \tilde{\tau}]\times[0,1])\cap C^{0,1}_{loc} ((0, \tilde{\tau})\times(0,1))$ to a solution $\tilde{u}$ of \eqref{eq1}. This was actually proved for 
$u_{\eps,0}\equiv u_0$, but an inspection of the proof shows that this is true in the general case. 
 The uniqueness of the solution of \eqref{eq1} implies that $\tilde{u}=u$ and that the whole sequence converges to $u$.
  We recall that $u_{\eps}$ is bounded in $L^{\infty}\left([0, \tilde{\tau}]; W^{1, \infty}(0, 1)\right)$  (see Step 3 of the proof of  Theorem  1.1 \cite{amal}). The boundary regularity result of Proposition \ref{friedmanetendu} implies that the convergence of  $\left\{(u_{\eps})_x\right\}$ to $u_x$ holds in $C_{loc}((0,\tilde{\tau})\times[0,1])$ (that is up to the boundary). Now fix $T\in (0, T_{max})$ and let 
\begin{eqnarray*}\tilde{T}&:=&
\text{sup}\left\{ s>0 \quad\text{such that}\quad T(u_{\eps})>s \quad\text{for}\quad \eps >0 \quad \text{small and}\right.\\
&&\left.\qquad u_{\eps}\to u \quad\text{in}\quad C([0, s]\times[0,1])\cap C^{0,1}_{loc} ((0, s)\times[0,1])\right\}.
\end{eqnarray*}
We know that $\tilde{T}\geq \tilde{\tau}(A)>0$. Assume that $\tilde{T}<T$. 
Set $A_1=\underset{t\in [0,T]}{\sup}\, \|u(t)\|_{W^{1,\infty}}<\infty$. For any $\eta\in (0,\tilde {T})$, we have 
\begin{equation}\label{enfind}
u_\eps(\tilde{T}-\eta) \to u(\tilde{T}-\eta) \qquad\text{in}\quad W^{1,\infty}(0,1).
\end{equation}
Thanks to \eqref{enfind} and the small-time existence and convergence result mentioned at the beginning of the proof,  we can find $\tau=\tau(A_1)>0$ (independent of $\eta$) and $\eps_0=\eps_0(\eta)>0$ such that the problem 
\begin{eqnarray}
\begin{array}{lll}
(u^{\eta}_{\varepsilon})_t=\big( \left(|(u^{\eta}_{\varepsilon})_x|^2 +\varepsilon^2 \right)^{\frac{p-2}{2}}(u^{\eta}_{\varepsilon})_{x}\big)_x+B_{\eps}((u^{\eta}_{\varepsilon})_x)
 & t>0, & x\in (0,1),\\
u^{\eta}_{\varepsilon}(t,0)=0,  \qquad u^{\eta}_{\varepsilon}(t,1)=M,     &t>0,& \\
u^{\eta}_{\varepsilon}(0,x)=u_{\eps}(\tilde{T}-\eta, x), &  & x\in(0,1),
       \end{array}
 \end{eqnarray}
admits  a unique classical solution $u^{\eta}_{\eps}$ on $[0, \tau]$. 
Moreover, we have $u_{\eps}^{\eta}\to u(\tilde{T}-\eta+\cdot,\cdot)$ in $C([0,\tau]\times[0,1])\cap C^{0, 1}((0,\tau)\times[0,1])$.
We can extend the solution $u_{\eps}$ of \eqref{eq1a} on $[0, \tilde{T}-\eta+\tau]$ by setting $u_{\eps}(t,x)=\left\{\begin{array}{ll}
u_{\eps}(t,x)\quad \text{for} \quad x\in[0, \tilde{T}-\eta],\\
u^{\eta}_{\eps}(t,x)\quad \text{for} \quad x\in [\tilde{T}-\eta, \tilde{T}-\eta+\tau]\end{array}\right.$.

It follows that $u_{\eps}\to u$ in $C([0,\tilde{T}-\eta+\tau]\times[0,1])\cap C^{0, 1}((0,\tilde{T}-\eta+\tau)\times[0,1])$.
Since $\tilde{T}-\eta+\tau>\tilde{T}$ for $\eta$ small enough, this contradicts the definition of $\tilde{T}$.
The second assertion follows from the estimates given in \cite[Inequalities 2.16 and 2.19]{amal}.
\end{proof}

Let us also note that due to $q> p> 2$, we have for $\varepsilon$ small enough  $$(p-1) \varepsilon^p \cosh(\varepsilon x)^{p-1}\ge \varepsilon^q \cosh(\varepsilon x)^q$$ (it suffices to take $0<\varepsilon <\cosh(1)^{\frac{p-1-q}{q-p}}$).  Hence  $ \left\| u_0\right\|_{L^{\infty}}+M+2-\cosh( \varepsilon x)$ is a supersolution for problem \eqref{eq1a}. It is also easy to see that $-\left\|u_0\right\|_{L^{\infty}}$ is a subsolution. Therefore there exists $K>0$ depending only on $\left\|u_0\right\|_{L^{\infty}}$ such  that, \begin{equation}\label{ahphil}
\forall \varepsilon\in (0,1/2), \qquad\qquad \left\| u_{\varepsilon}(t,x)\right\|_{L^{\infty}}\leq K.
\end{equation}
\subsection{Construction of the Lyapunov functional}
Now we construct a Lyapunov functional for \eqref{eq1a}  with the help of the
technique developed by Zelenyak \cite{zelen}. Let $D_K=[-K, K]\times \mathbb{R}$, where $K$ is  the constant in \eqref{ahphil}. We look for a pair of functions $\Phi_{\varepsilon}\in C^1(D_{K};\mathbb{R})$ and  $\Psi_{\varepsilon}\in C(D_K; (0,\infty))$ with the following property:\\
For any solution $u_{\varepsilon}$ of \eqref{eq1a} with $|u_{\varepsilon}|\leq K$, defining\\

$$\mathcal{L}_{\varepsilon}(u_{\epsilon}(t))=\int_0^1 \Phi_{\varepsilon} \left( u_{\varepsilon}(t,x), (u_{\varepsilon})_x(t,x)\right) dx,$$
it holds
$$\dfrac{d}{dt}\mathcal{L}_{\varepsilon}(u_{\varepsilon}(t))= -\int_0^1 \Psi_{\varepsilon}\left( u_{\varepsilon}(t,x),(u_{\varepsilon})_x(t,x) \right) \, (u_{\varepsilon})_t^2(t,x) \,dx.$$
Since $(u_{\varepsilon})_t(t,0)=(u_{\varepsilon})_t(t,1)=0$, we have

\begin{align*}
&\dfrac{d}{dt}\int^1_0 \Phi_{\varepsilon}(u_{\varepsilon}, (u_{\varepsilon})_x) dx=\int_0^1  (u_{\varepsilon})_t\cdot  (\Phi_{\varepsilon})_u \left(u_{\varepsilon}, (u_{\varepsilon})_x\right)+ (u_{\varepsilon})_{x t}\cdot  (\Phi_{\varepsilon})_v\left( u_{\varepsilon}, (u_{\varepsilon})_x\right) dx\\
&\quad=\int^1_0 (u_{\varepsilon})_t\big[(\Phi_{\varepsilon})_u\left(u_{\varepsilon}, (u_{\varepsilon})_x\right) 
 - (u_{\varepsilon})_x \cdot(\Phi_{\varepsilon})_{u v}\left(u_{\varepsilon}, (u_{\varepsilon})_x\right)-(u_{\varepsilon})_{x x}\cdot(\Phi_{\varepsilon})_{vv}\left(u_{\varepsilon}, (u_{\varepsilon})_x\right)\big] dx.
\end{align*}
So it is natural to require that
\begin{align*}
&(\Phi_{\varepsilon})_u\left(u_{\varepsilon}, (u_{\varepsilon})_x\right) 
- (u_{\varepsilon})_x\cdot (\Phi_{\varepsilon})_{uv}\left(u_{\varepsilon}, (u_{\varepsilon})_x\right)-(u_{\varepsilon})_{xx}\cdot(\Phi_{\varepsilon})_{vv}\left(u_{\varepsilon}, (u_{\varepsilon})_x\right)\\
&\qquad\qquad=-\Psi_{\varepsilon}(u_{\varepsilon},(u_{\varepsilon})_x) \cdot(u_{\varepsilon})_t\\
&\qquad\qquad=-\Psi_{\varepsilon}(u_{\varepsilon},(u_{\varepsilon})_x)\left[(p-1) \left(   |(u_{\varepsilon})_x|^2 +\varepsilon^2 \right)^{\frac{p-4}{2}}\left( |(u_{\varepsilon})_x|^2 +\frac{\varepsilon^2}{p-1} \right)(u_{\varepsilon})_{xx}\right.\\
&\qquad\qquad+\left.\left(|(u_{\varepsilon})_x|^2+\varepsilon^2\right)^{\frac{q-2}{2}}
\left(|(u_{\varepsilon})_x|^2 +
\frac{\varepsilon^2}{p-1} \right)\right]
\end{align*}
A sufficient condition is
\begin{align}
&(\Phi_{\varepsilon})_{vv}(u,v)=(p-1)\Psi_{\varepsilon}(u,v) \left(v^2+\varepsilon^2\right)^{\frac{p-4}{2}}\left(v^2+\frac{\eps^2}{p-1}\right),\label{condi1}\\
&(\Phi_{\varepsilon})_u(u,v)- v(\Phi_{\varepsilon})_{uv}(u,v)=-\Psi_{\varepsilon}(u,v) \left(v^2+\varepsilon^2\right)^{\frac{q-2}{2}}\left(v^2+\frac{\eps^2}{p-1}\right),\label{condi2}
\end{align}
 that is $\Phi_{\varepsilon}$ satisfies the differential equation:
\begin{equation}\label{eqdif}
(\Phi_{\varepsilon})_u(u,v)- v(\Phi_{\varepsilon})_{uv}(u,v)+\dfrac{\left(v^2+\varepsilon^2 \right)^{\frac{q-p+2}{2}}}{p-1}(\Phi_{\varepsilon})_{vv}(u,v)=0.
\end{equation}
We follow the method used in \cite{Arri} to find such nice functions.
For a given function $\rho_{\varepsilon}(u,v)$, let us denote
$$ H_{\varepsilon}=(\rho_{\varepsilon})_u+\dfrac{\left(v^2+\varepsilon^2 \right)^{\frac{q-p+2}{2}}}{p-1}(\rho_{\varepsilon})_{vv}-v (\rho_{\varepsilon})_{uv}.$$
Here we assume that $\rho_{\varepsilon}, (\rho_{\varepsilon})_u, (\rho_{\varepsilon})_v ,(\rho_{\varepsilon})_{uv}$ are continuous and  $C^1$  in $v $ in $D_K$, and that
$(\rho_{\varepsilon})_{vv}$ is continuous in $D_K$ and, except perhaps at $v = 0$, $C^1$ in $v$.\\
We want to have  $( H_{\varepsilon})_v=0$,  so that  $H_{\varepsilon}(u,v)=H_{\varepsilon}(u,0)=H_{\varepsilon}(u)$.
We compute
$$ ( H_{\varepsilon})_v=\dfrac{ \left( v^2+\varepsilon^2 \right)^{\frac{q-p+2}{2}}}{p-1}(\rho_{\varepsilon})_{vvv}+\left(\dfrac{q-p+2}{p-1}\right) v \left(v^2+\varepsilon^2\right)^{\frac{q-p}{2}} (\rho_{\varepsilon})_{vv}-v (\rho_{\varepsilon})_{uvv}.$$
To have $ ( H_{\varepsilon})_v= 0$, it suffices that $f_{\varepsilon}=(\rho_{\varepsilon})_{vv} $ satisfies the following conditions:
\begin{eqnarray} \label{tou}
\left\{
\begin{array}{ll}
   (f_{\varepsilon})_u-\left(\dfrac{q-p+2}{p-1}\right) \left(v^2+\varepsilon^2 \right)^{\frac{q-p}{2}} f_{\varepsilon} -
\dfrac{\left(v^2+\varepsilon^2\right)^{\frac{q-p+2}{2}}}{(p-1)v} (f_{\varepsilon})_v=0 \quad |u|\leq K ,  v\neq 0,\\
     (f_{\varepsilon})_v(u,0)=0.
   \end{array}
              \right.
\end{eqnarray}
Now, the equation (\ref{tou})  can be solved by the method of characteristics.
For any $K>0$ such that $|u|\leq K$, one finds that the function  defined by
$$f_{\varepsilon}(u,v)=\left[1+\left(\frac{q-p}{p-1}\right)\left(v^2+\varepsilon^2\right)^{\frac{q-p}{2}}\left(K+1-u\right)\right]^{-\frac{q-p+2}{q-p}}>0 $$
is a solution of  (\ref{tou}) on $[-K,K]\times \mathbb{R}$.
Define $\rho_{\varepsilon}$ by
$$\rho_{\varepsilon}(u,v)=\int_0^v \int_0^z f_{\varepsilon}(u,s)\  ds\  dz \geq 0,$$
and let then
 \begin{equation}
 \Phi_{\varepsilon} (u,v)=\rho_{\varepsilon}(u,v)-\int_0^u H_{\varepsilon}(s,0) \ ds + K+1.
 \end{equation}
We added the constant $K+1$ to ensure that $\Phi_{\varepsilon} \geq 0$. In fact, given that $\varepsilon\leq 1/2$, $2<p$ and  $0\leq(\rho_{\varepsilon})_{vv}\leq 1$, we get
  $(\rho_{\varepsilon})_u(s,0)=0$ and
\begin{equation}
0\leq H_{\varepsilon}(s,0)=(\rho_{\varepsilon})_u(s,0)+\dfrac{\varepsilon^{q-p+2}}{p-1} f_{\varepsilon}(s, 0)\leq 1.
\end{equation}
Consequently, using that $|u|\leq K$, we get
\begin{equation*}
\begin{split}
 -\int_0^u H_{\varepsilon}(s,0) \ ds\geq -u \,\geq -K\quad &\mathrm{for} \quad u \in [0, K],\\
-\int_0^u H_{\varepsilon}(s,0) \ ds\geq 0\quad &\mathrm{for} \quad u \in [-K, 0].
\end{split}
\end{equation*}
Using the definition of $H_{\varepsilon}$ and the fact that $H_{\varepsilon}(u,v)=H_{\varepsilon}(u,0)$, we see that:
$$(\Phi_{\varepsilon})_u- v(\Phi_{\varepsilon})_{uv}(u,v)+\dfrac{\left(v^2+\varepsilon^2 \right)^{\frac{q-p+2}{2}}}{p-1}(\Phi_{\varepsilon})_{vv}(u,v)=0,$$
i.e. $\Phi_{\varepsilon}$ satisfies \eqref{eqdif}, hence \eqref{condi1}-\eqref{condi2} with
 \begin{equation}
\Psi_{\varepsilon}(u,v)=\left( v^2+\dfrac{\eps^2}{p-1}\right)^{-1}\dfrac{\left(v^2+\varepsilon^2\right)^{\frac{4-p}{2}}(\rho_{\varepsilon})_{vv}}{p-1}\geq\dfrac{\left(v^2+\varepsilon^2\right)^{\frac{2-p}{2}}(\rho_{\varepsilon})_{vv}}{p-1}>0.
\end{equation}
It follows that
$$\dfrac{d}{dt}\mathcal{L}_{\varepsilon}(u_{\varepsilon}(t))
 =-\int_0^1 \dfrac{ \left((u_{\varepsilon})_x^2+\varepsilon^2\right)^{\frac{2-p}{2}}(\rho_{\varepsilon})_{vv}}{(p-1) \left((u_{\varepsilon})_x^2+\dfrac{\varepsilon^2}{p-1}\right)}\, (u_{\varepsilon})_t^2\  dx=-\int_0^1\Psi_{\varepsilon}(u_{\varepsilon}, (u_{\varepsilon})_x) \, (u_{\varepsilon})_t^2 \ dx.$$
Due to $q>p>2$, we remark that, $\forall \eps\in(0,1)$, $|u|\leq K$ and $v\in\mathbb{R}$,
\begin{eqnarray}\label{pousha}
 \Psi_{\varepsilon}(u,v)\geq \mathcal{A}(v)=\dfrac{\left(v^2+1\right)^{\frac{2-p}{2}}}{p-1} \left[1+\frac{(q-p)}{(p-1)}\left(v^2+1\right)^{\frac{q-p}{2}}\left(2K+1\right)\right]^{-\frac{q-p+2}{q-p}}.
\end{eqnarray}
As a  consequence of the existence of the approximate Lyapunov functional, we have the following estimate.
\begin{prop}\label{estimationfonc}
Assume that  $q>p>2$ and let  $u$ be a global weak solution  of \eqref{eq1}. Then for any $T>1$ and $\delta>0$,  There exists $ C=C(\left\|u_0\right\|_{W^{1,\infty}}, \delta,p,q)>0$ such that 
\begin{equation}\label{seb}
 \int_1^T\int_{\delta}^{1-\delta} (u_t)^2 dx dt\leq C.
\end{equation}
\end{prop}

\begin{proof}

First let us remark that Lemma \ref{as}  implies that, for any $\delta>0$,
 
\begin{equation}\label{refe}
 | u_x| \leq  C(\delta) \qquad \mathrm{in}\quad [1, \infty)\times[\delta, 1-\delta].
\end{equation}
Now we fix $T>1$ and $\delta\in(0, 1/2)$. 
On the one hand, by \eqref{pousha}, we have
\begin{eqnarray}
  \int_0^T\int_{\delta}^{1-\delta}\mathcal{A}\left((u_{\varepsilon})_x\right)\cdot (u_{\varepsilon})_t^2(t,x)\, dx\, dt
&\leq&  \int_{0}^{T}\int_{0}^1\Psi_{\varepsilon}(u_{\varepsilon}, (u_{\varepsilon})_x)\cdot(u_{\varepsilon})_t^2 (x,t)\,dx \,dt \nonumber\\
&=& \mathcal{L}_{\varepsilon}(u(0))-\mathcal{L}_{\varepsilon}(u_{\varepsilon}(T))\label{fifou}\\
&\leq& \tilde{C}(\left\|u_0\right\|_{W^{1,\infty}}).\nonumber
\end{eqnarray}

On the other hand, by Proposition \ref{estimsuite}, there exists $\eps_0(\delta, T)$ such that, for all $\eps<\eps_0$, $x\in[0,1]$, $t\in[1,T]$,
$$|(u_{\eps})_x(t,x)-u_x(t,x)|\leq C(\delta).$$
Then, by \eqref{refe}, $|(u_{\eps})_x|\leq 2C(\delta)$ for $(t,x)\in[1,T]\times[\delta, 1-\delta]$ so that

\begin{eqnarray*}
 \int_0^T\int_{\delta}^{1-\delta}\mathcal{A}\left((u_{\varepsilon})_x\right)\cdot (u_{\varepsilon})_t^2(t,x)\, dx\, dt &\geq &\int_1^T\int_{\delta}^{1-\delta}\mathcal{A}\left((u_{\varepsilon})_x\right)\cdot (u_{\varepsilon})_t^2(t,x)\, dx\, dt\\
&\geq &\theta(2C(\delta)) \int_{1}^{T}\int_{\delta}^{1-\delta}(u_{\varepsilon})_t^2 (x,t)\,dx \,dt,
\end{eqnarray*}
where $\theta(R)=\text{inf} \left\{ \mathcal{A}(v); \quad |v|\leq R\right\}>0$.
Letting $\eps\to 0$ and using a lower semicontinuity argument as well as \eqref{fifou}, we obtain
\begin{equation}\label{hihi}
\theta(2C(\delta))\int_{1}^{T}\int_{\delta}^{1-\delta} (u)_t^2(t,x)\, dx\, dt
\leq \tilde{C}(\left\|u_0\right\|_{W^{1,\infty}}),
\end{equation}
The result immediately follows.
\end{proof}
\subsection{Convergence to steady states}
\begin{prop}\label{converg}
Let $u$ be a global weak solution of (\ref{eq1}). Then $M\leq M_b$ and  $u(t)$ converges in $C([0,1])$ to a steady state of \eqref{eq1} as $t\rightarrow \infty$. Moreover the convergence also holds in $C^{1}([\delta, 1-\delta])$ for all $\delta>0$.
\end{prop}
\begin{Proof}
 Assume that $u$ is a global weak solution of \eqref{eq1}. Fix a sequence  $(t_k)_{k \in \mathbb{N}}$, $1\leq t_k\to\infty$ and set $w_k (t,x)=u(t+ t_k,x)$.
By \eqref{compaprin}, we know that 
\begin{equation}\label{not1}
 |u|\leq \max\left\{\left\|u_0\right\|_{\infty}, M\right\}\qquad \quad\mathrm{in}\quad  [1,\infty)\times [0,1],
\end{equation}
Using lemma \ref{as}  we have 
\begin{equation}\label{referee}
 | u_x|\leq  C(\delta), \qquad \mathrm{in}\quad [1, \infty)\times[\delta/2, 1-\delta/2].
\end{equation}
Thus applying a result of DiBenedetto-Friedman \cite{dib}, we have that $\left\{w_k\right\}$ and $\left\{(w_k)_x\right\}$ are H\"older continuous in $[\delta,T-\delta]\times[\delta, 1-\delta]$ with a H\"older norm independent of $k$.
It follows that $\left\{ w_k \right\} $ and $\left\{ (w_k)_x \right\} $ are relatively compact in $C\left([\delta,T-\delta]\times[\delta, 1-\delta]\right)$ for any $\delta$, $T>0$.
Thus, by the Arzel\`a-Ascoli theorem and a diagonal procedure, there exist a subsequence $(t_{k_l})_{l \in\mathbb{N}}$ of $(t_k)$  and a function $W\in C\left((0, \infty)\times(0, 1)\right)$, $ W_x\in C\left((0,\infty)\times(0,1)\right)$ such that for any $\delta, T>0$
\begin{align}
 w_{k_l}\to W \qquad  \text{strongly  in} \quad  C\left([\delta, T-\delta]\times[\delta,1-\delta]\right) \qquad  &\mathrm{as}\quad  l\to\infty.\\
 ( w_{k_l})_x\to W_x \qquad \text{strongly in} \quad  C\left([\delta,T-\delta]\times[\delta, 1-\delta]\right)\qquad  &\mathrm{as} \quad l \to\infty.
\end{align}
and $W$  is a distributional solution of
\begin{equation*}
W_t-(|W_x|^{p-2} W_x)_x=|W_x|^{q},\qquad    t>0,    x\in (0, 1).
 \end{equation*}
Further, using lemma \ref{as} and $ q> p$, we get that for some $r>1$
\begin{equation}\label{not2}
 \left\|(w_k)_x\right\|_{L^{\infty}(1,\infty; L^r(0,1))}\leq C.
\end{equation} 
Combining \eqref{not1} with \eqref{not2}, we get that, for  each fixed $t>0$, $w_k(t,\cdot)$ is relatively compact in $C([0, 1])$. Consequently for any $t>0$, $W(t, .)$ can be extended to a continuous function on $[0, 1]$ satisfying
\begin{equation*}
 W(t, 0)=0\qquad \qquad W(t,1)= M.
\end{equation*}
Proposition \ref{estimationfonc} implies that \[\displaystyle \int_0^{\infty}\int_{\delta}^{1-\delta} ( w_ {k_l})_t^2 (t,x) \ dx\, dt\to \, 0, \qquad\mathrm{as}\quad l\to\infty.\]
Since $( w_{k_l})_t\to W_t$ in $\mathcal{D}'( (0, \infty)\times(0, 1))$ and $\delta\in(0,1)$ is arbitrary, it follows that $ W_t\equiv 0$. Thus $W$ is a  steady state of (\ref{eq1}) which implies that $M\leq M_b$. Given that the sequence $t_{k}\to\infty$ is arbitrary and the steady states (for given $M$) are unique, it follows that the whole solution $u(t)$ converges to $W$.

\end{Proof}

\section{Proof of Theorem \ref{Main result}}

\subsection{GBU  profiles and lower bound on  $u_x$}

Thanks to \eqref{derivtinf} in   Theorem \ref{reg}, we shall derive the following lemma  providing a lower bound on the blow up profile of $u_x$ in case  GBU occurs  in finite or infinite  time near $x=0$ or $1$.

\begin{lem}\label{mariep} Let  $u$ be a global  weak solution of \eqref{eq1} and $t_0>0$. Let $C_1>0$ be the constant given in the estimate \eqref{estimut} of Theorem \ref{reg}. There exist $ C_4=C_4( C_1, p, q)$, $ C_5=C_5(p,q) >0$ with the following property. For all   $t\in [t_0,+\infty)$ and  $0\leq y\leq  x\leq 1$
\begin{equation}\label{df}
 \big[ |u_x|^{p-2}u_x^+ (t,x)+C_4 \big]^{\frac{p-1-q}{p-1}}\leq \big[|u_x|^{p-2}u_x^+(t,y)+C_4\big]^{\frac{p-1-q}{p-1}} +C_5 (x-y),
\end{equation}
and 
\begin{equation}\label{deb}
 \big[ |u_x|^{p-2}(-u_x)^+ (t,1-x)+C_4 \big]^{\frac{p-1-q}{p-1}}\leq \big[|u_x|^{p-2}(-u_x)^+(t,1-y))+C_4\big]^{\frac{p-1-q}{p-1}} +C_5 (x-y).
\end{equation}
\end{lem}
\begin{proof} Fix $t \in [t_0, T_{max})$. and let  $z(x)=|u_x |^{p-2}u_x^+(t,x) +C_1^{\frac{p-1}{q}}$, where  $C_1$ is given by the estimate of  $|u_t|$ in Theorem \ref{reg}.
Using  that $|u_t|\leq  C_1$ and $|u_x(t)|^q$ is bounded,   we get that, $z$ is Lipschitz in $(0,1)$ with  a Lipschitz  bound  depending on $t, p,q,  ||u_0||_{\infty}$ and $||u_x(t)||_{\infty}^q$. Moreover  the  function  $z$ satisfies almost everywhere
$$\begin{array}{lll}
z'+z^{q/(p-1)}&=\left(|u_x|^{p-2} u_x(t,x)\right)_x\mathbf{1}_{\left\{u_x>0\right\}}+\left(|u_x|^{p-2}u_x^+(t,x)+C_1 ^{\frac{p-1}{q}}\right)^{\frac{q}{p-1}} \\
& \geq\left[(|u_x|^{p-2}u_x(t,x))_x+|u_x|^q \right]{\mathbf 1}_{\left\{u_x>0\right\}}+C_1\\
&\geq 0.
\end{array}$$
For $0\leq y\leq x\leq 1$, an integration yields  $$z(x)^{(p-1-q)/(p-1)}\leq z(y)^{(p-1-q)/(p-1)}+\left(\dfrac{q-p+1}{p-1}\right)(x-y),$$
 that is \eqref{df} with $C_4=C_1^{\frac{p-1}{q}}$  and $C_5=\dfrac{q-p+1}{p-1}$.

 The estimate \eqref{deb} can be obtained similarly by considering  $z(x)=|u_x|^{p-2}(-u_x)^+ (t,1-x)+C_1^{\frac{p-1}{q}}$.
\end{proof}
\begin{rem}
 Lemma \ref{mariep} yields in particular a lower bound on the gradient blow-up profile,
which complements the upper bounds in \eqref{fer}-\eqref{erf}. 
Namely, if $x=0$ is a GBU point (in finite or infinite time), i.e.
if $|u_x|$ is unbounded in any neighborhood of $T_{max}$ and $0$, then 
$$\limsup_{t\to T_{max}} u_x(t,x)\geq C(p,q)x^{-1/(q-p+1)}$$
for all sufficiently small $x>0$.
The analogous estimate holds if $x=1$ is a GBU point.
\end{rem}
Now let us state the following lemma which is a direct consequence of the convergence of $u$ to the steady state.
\begin{lem}\label{cantine}
Let $M\geq 0$ and let $u$  be a global weak solution of \eqref{eq1}.Then it holds 
\begin{equation}
\lim_{t\to +\infty} \left(\underset{[0,1]}{\max} \, u(t,x) \right)=M.
\end{equation}
\end{lem}
\begin{proof}
Since $u(t,1)=M$, we get $\left(\underset{[0,1]}{\max} \, u(t,x) \right)\geq M$. Next, using that $u(t)\to W$ in $C([0,1])$  (see Proposition \ref{converg}) and $W\leq M$, it holds that
$\forall \eps>0, \exists t_{\eps}>0$ such  that if $t>t_{\eps}$ then
$$u(t,x)\leq W(x)+\eps\leq M+\eps\qquad x\in[0,1].$$
It results that $\underset{[0,1]}{\max}\, u(t,x)\leq M+\eps$ if $t>t_{\eps}$.
\end{proof}

Thanks to this property we can rule out infinite time gradient blow-up towards $-\infty$ when $x\to 1$.
\begin{lem}\label{casfacil}
Let $u$ be a global weak solution of \eqref{eq1}. Then

\begin{equation}
\underset{[0, \infty)\times (0,1)}{\inf}\, u_x>-\infty.
\end{equation}

\end{lem}
\begin{proof}
We proceed by contradiction. Assume that the lemma is false. Then, by Lemma \ref{as},
there exist a sequence $t_n\to +\infty$ and $x_n\to 0$ such that $u_x(t_n, 1-x_n)\to -\infty$.
Fix $\eps>0$, then for $n\geq n(\eps)$ large enough, we have $x_n<\varepsilon$ and
$$|u_x|^{p-2}(-u_x)^+(t_n, 1-x_n)\geq \eps^{-(p-1)/(q-p+1)}.$$
Taking $t=t_n$ and $y=x_n$ in  \eqref{deb}, we get that for  $n\geq n_0(\eps)$ large enough, we have  for $x_n\leq x\leq \eps$
\begin{eqnarray*}
\big[ |u_x|^{p-2}(-u_x)^+ (t_n,1-x)+C_4 \big]^{\frac{p-1-q}{p-1}}&\leq& \big[|u_x|^{p-2}(-u_x)^+(t_n,1-x_n)+C_4\big]^{\frac{p-1-q}{p-1}} +C_5x\\
&\leq&  (C_5+1)\eps.
\end{eqnarray*}
This implies that
\begin{equation}
\big[ |u_x|^{p-2}(-u_x)^+ (t_n,1-x) \geq \left((C_5+1)\eps\right)^{(1-p)/(q-p+1)}-C_4, \qquad x_n\leq x\leq \eps.
\end{equation} 
Choosing $\eps=\eps (C_5, C_4)>0$ small enough,  we get that $u_x(t_n,1-x)\leq -1$ on $[x_n, \eps]$, hence 
$$u(t_n, 1-x)\geq u(t_n, 1-x_n)+ (x-x_n), \qquad x_n\leq x\leq \eps.$$
Using that $u(t_n, 1-x_n)\to M$ (by Proposition \ref{converg}) and recalling  Lemma \ref{cantine},  we end up with a contradiction.
\end{proof}

 \begin{rem}
Thanks to lemma \ref{casfacil} we deduce that, for the case $M=M_b$, if there exist global  solutions with infinite time gradient blow up  (we expect that  this could occur for some particular initial data), then $u_x$ can only blow up at $x=0$.
\end{rem}

\subsection{Completion of the proof of Theorem \ref{Main result}}\label{secconv}
\subsection*{Proof of the boundedness of $u_x$ for $M=0$}
Lemma \ref{casfacil} is sufficient to prove the main theorem in the case $M=0$.
Let $u$ be a global solution of \eqref{eq1}. For $M=0$, we note that  $w(t,x):=u(t, 1-x)$ solves \eqref{eq1} with $u_0(1-x)$ as initial data. Lemma \ref{casfacil} implies that $u_x$ and $w_x$ are bounded below on $[0, +\infty)\times (0,1)$, therefore $u_x$ is bounded.  See the Subsection \ref{secconv} for the proof of the convergence to the steady state in the $W^{1, \infty}(0, 1)$ norm.
 \hfill$\square$

\subsubsection*{Proof of the boundedness of $u_x$ for $0< M<M_b$}
We proceed by contradiction. Assume that $u$ is a global weak solution which is unbounded in  $W^{1,\infty}$. We know that when $t\rightarrow\infty$, $u$ converges to $W=W_M$ in $C[0,1]$ and in $C^1[\delta,1-\delta]$ for all $\delta>0$.
Since  $u_x$ is unbounded  and can only blow up to $+\infty$ at $x=0$, there exist  sequences $t_n\rightarrow\infty$, $x_n\to 0$ such that
\begin{equation}
u_x(t_n,x_n)\rightarrow + \infty
\end{equation}
Taking $t=t_n$ and $y=x_n$ in  \eqref{df} and sending $n\to\infty$, we deduce that, for any $x\in (0, 1)$
\begin{equation*}
\big[|W_x |^{p-2}W_x +C_4 \big]^{\frac{p-1-q}{p-1}}\leq C_5 \,x.
\end{equation*}
This would imply that
\begin{equation*}
 |W_x|^{p-2} W_x +C_4\geq \left(C_5 \, x \right)^{\frac{1-p}{q-p+1}}.
\end{equation*}
Passing to the limit  $x\rightarrow 0$  we get a contradiction since  $W=W_M\in C^1([0,1])$.
So all the global solutions are bounded in $W^{1,\infty}$.

\subsubsection*{Proof of the convergence in $C^1$ norm for $M\in [0, M_b)$}
This follows from the proof of Proposition \ref{converg}, with \eqref{refe} replaced by the boundedness of $u_x$ on  $[0, \infty)\times[0,1]$,
and using Proposition \ref{friedmanetendu} which is an extension of the result in  \cite{dib}.

\subsubsection*{Proof of Theorem \ref{Main result} for $M>M_b$}
This is an immediate consequence of Proposition \ref{converg}  and the fact that (\ref{stat}) admits no solution for $M > M_b$.\hfill$\square$
\subsubsection*{Further regularity for global solutions  for $0<M<M_b$}
As a consequence of the convergence of global solutions to the steady state  in $C^1[(0, 1)]$,  we have  the following proposition which is of independent interest. It  gives a result of further regularity of global solutions for large time. It is unknown whether or not such property is true in the case $M=0$.
\begin{prop}\label{miloun}
 Assume  that $0<M<M_b$ and let $u$ be a global weak solution of \eqref{eq1}.
Then there exist $\tilde{T} >0$ and $\tilde{\eta} >0$ such that
$$u_x\geq \tilde{\eta}	\qquad\text{on}\quad	[\tilde{T}, +\infty) \times [0, 1].$$
Moreover, $u$ becomes a classical solution on $[\tilde{T}, +\infty) \times [0, 1]$ 
\end{prop}

\begin{proof}
First let us note that there exists  $\eta >0$ such that $(W_{M})_x\geq 2 \eta > 0$ in $[0,1]$. Next,  by Theorem \ref{Main result}, we know that  $u_x\to W_x$ uniformly on $[0,1]$. Hence,  there exists $\tilde{T}>0$ such that
\begin{equation}\label{phil1}
u_x(t,x)> (W_M)_x(x)-\eta \geq\eta \qquad  \text{for all}\quad x\in [0, 1],\,  t> \tilde{T}.
\end{equation}
The last inequality implies that the differential equation is uniformly parabolic for $(t,x)\in [\tilde{T}, \infty]\times [0, 1]$. Hence, by the standard theory (see \cite{lady}) we know that $u\in C^{1,2}((\tilde{T}, \infty)\times[0,1])$
\end{proof}

\section*{Appendix}

\subsection*{Proof of Proposition \ref{friedmanetendu} on the regularity of the derivative up to the boundary}
Let $\eps\in[0,1), u_0\in W^{1,\infty},C>0$ and $F\in L^{r} ((0,T)\times(0,1))$ for some $r>2$ with $\left\|F\right\|_{L^{r}((0,T)\times(0,1))}\leq C$.  
Since $v_t\in L^2((0,T)\times (0,1))$ by assumption and $r>2$, it follows that
\begin{equation}\label{francoisj}
\left((|v_x|^2+\eps^2)^{\frac{p-2}{2}}v_x\right)_x \in L^2((0,T)\times (0,1)),
\end{equation}
and that the partial differential equation in \eqref{friedequa} is satisfied
in the sense of equality of functions in $L^2((0,T)\times (0,1))$.
 
Next, we define an extension $v^*$ of $v$ to $[-1, 2]$ by setting
\begin{equation}
v^*(t,x)=\left\{
\begin{array}{lll}
&-v(t,-x)&\qquad \text{if}\quad x\in [-1, 0)\\
&v(t,x)&\qquad \text{if}\quad x\in [0, 1]\\
&2M-v(t,2-x)&\qquad \text{if}\quad x\in (1, 2]
\end{array}
\right.
\end{equation}
We denote by $u_0^*$ the extension of $u_0$ to $[-1, 2]$. 
We will prove that $v^*$ is a weak solution of the following problem
\begin{eqnarray}\label{eq1aetendue}
\left\{
\begin{array}{lll}
 v^*_t=\left(\left(|v^*_x|^2 +\varepsilon^2 \right)^{\frac{p-2}{2}}v^*_{x}\right)_x +\tilde{F}(t,x)
 & t>0, & x\in (-1,2),\\
v^*(t,-1)=-M,  \qquad v^*(t,2)=2M,     &t>0,& \\
v^*(0,x)=u^*_{0}(x), &  & x\in(-1,2).
       \end{array}
\right.
 \end{eqnarray}
where $$\tilde{F}(t,x)=\left\{\begin{array}{lll}
&-F(t,-x)&\qquad \text{if}\quad x\in [-1, 0)\\
&+F(t,x)&\qquad \text{if}\quad x\in [0, 1]\\
&-F(t,2-x)&\qquad \text{if}\quad x\in (1, 2]
\end{array}
\right.$$
Indeed, let $\psi\in C([0, \tau]\times[-1,2]) \cap L^p((0,\tau); W^{1,p}_0((-1,2))$. 
Due to \eqref{francoisj}, for a.e. $t\in (0,T)$, we have $(|v_x|^2+\eps^2)^{\frac{p-2}{2}}v_x(t,\cdot) \in W^{1,2}(0,1)\subset C([0,1])$. By elementary distribution theory (jump formula), it readily follows that $\left(|v^*_x|^2+\eps^2\right)^{\frac{p-2}{2}}v^*_x(t,\cdot) \in W^{1,2}(-1,2)\subset C([-1,2])$.
For a.e. $t\in (0,T)$, we can thus write:

\begin{eqnarray*}
 \int_0^1 v^*_t(t,x) \psi(t,x) \,dx &=&\int_0^1\left(\left(|v^*_x|^2 +\varepsilon^2 \right)^{\frac{p-2}{2}}v^*_{x}\right)_x(t,x) \psi(t,x) \, dx + \int _0^1 F(t,x)\,\psi(t,x) \,dx\\
&=&\left(|v^*_x|^2 +\varepsilon^2 \right)^{\frac{p-2}{2}}v^*_{x}(t,1) \psi(t,1)
-\left(|v^*_x|^2 +\varepsilon^2 \right)^{\frac{p-2}{2}}v^*_{x}(t,0)\psi(t,0)\\
&-&\int_0^1\left(|v^*_x|^2 +\varepsilon^2 \right)^{\frac{p-2}{2}}v^*_{x}(t,x)\,\psi_x(t,x) \,dx + \int _0^1 F(t,x)\,\psi(t,x) \,dx.
\end{eqnarray*}
Using that $\psi(t,-1)=0$ and $\psi(t,2)=0$, we have
\begin{eqnarray*}
 \int_{-1}^0 v^*_t(t,x)\, \psi(t,x)\, dx 
&=&\int_{-1}^0\left( \left(|v^*_x|^2 +\varepsilon^2 \right)^{\frac{p-2}{2}}v^*_{x}\right)_x(t,x)\, \psi(t,x) \, dx -\int _{-1}^0 F(t,-x)\,\psi \,dx\\
&=&\left(|v^*_x|^2 +\varepsilon^2 \right)^{\frac{p-2}{2}}v^*_{x}(t,0) \psi(t,0) - \int_{-1}^0F(t,-x)\psi(t,x)\, dx\\
&-&\int_{-1}^0\left(|v^*_x|^2 +\varepsilon^2 \right)^{\frac{p-2}{2}}v^*_{x}(t,x)\psi_x(t,x) \,dx,
\end{eqnarray*}
and 
\begin{equation*}
\begin{split}
\int_{1}^2 v^*_t(t,x)\psi(t,x) dx &=\int_{1}^2\left((|v^*_x|^2 +\varepsilon^2 )^{\frac{p-2}{2}}v^*_{x}\right)_x(t,x)\psi(t,x) dx-\int_1^2 F(t,2-x)\psi(t,x)dx\\
&=-\left(|v^*_x|^2 +\varepsilon^2 \right)^{\frac{p-2}{2}}v^*_{x}(t,1) \psi(t,1) -\int_{1}^2 F(t, 2-x)\,\psi(t,x) \,dx\\
&-\int_{1}^2\left(|v^*_x|^2 +\varepsilon^2 \right)^{\frac{p-2}{2}}v^*_{x}(t,x)\,\psi_x(t,x) dx.
\end{split}
\end{equation*}
Summing  these identities and integrating over $(0,T)$, it follows that
\begin{eqnarray}
 \int_0^{T}\int_{-1}^2 v^*_t\, \psi\, dx\,dt &+&
\int_0^{T}\int_{-1}^2\left(|v^*_x|^2 +\varepsilon^2 \right)^{\frac{p-2}{2}}v^*_{x}(t,x)\,\psi_x(t,x)\, dx\,dt\label{hellomilou}\\
&=&\int_0^{T}\int_{-1}^2 \tilde{F} (t,x)\,\psi \,dx\,dt.\nonumber
\end{eqnarray}
Next,   since $\left\|\tilde{F}(t,x)\right\|_{L^{r} ((0,T)\times(-1,2))}\leq C \left\|F\right\|_{L^{r}(0,T)\times(0,1))}$,  using  a result of DiBenedetto and Friedman (see \cite{dib} and \cite[chapter 9]{livdib}  for the case  $\eps>0$) on H\"older regularity of gradient of some degenerate parabolic problems, we get that, for any $\eta>0$, $v^*_x \in  C^{\alpha}_{loc}([\eta, T-\eta]\times(-1,2))$ where $\alpha>0$ and the norm of $v^*_x$ depend only on $\left\|F\right\|_{L^{r}(0,T)\times(0,1))}, \left\|v^*_x\right\|_{L^p}$ and  $\left\|v^*\right\|_{L^{\infty}_t,L^2_x}$. We get the desired result recalling that $v^*_x=v_x$ on $[0,1]$ and using that $[0,1]\subset (-1,2)$.

\bigskip
\textbf{Addendum to \cite{amal}}\\
In \cite{amal}, the  step 3  of the proof of the convergence of the approximate solutions relied on a  result on the H\"older regularity of the gradient of some degenerate parabolic systems. The author wants to point out that the correct reference should be \cite{livdib} where an extension of the results of \cite{dib} to the regularized problem is mentioned.
\bigskip

\textbf{Acknowledgments.} The author would like to thank Professor Ph. Souplet
for useful suggestions during the preparation of this paper. The author would also like to thank the referee for  carefully reading the manuscript and giving constructive comments which substantially helped improving the quality of the paper.
\nocite{*}
\bibliography{bibfinal}
\bibliographystyle{plain}
\end{document}